\documentclass[12pt]{amsart}
\usepackage[utf8]{inputenc}
\usepackage{amsmath}
\usepackage{amssymb}
\usepackage{amsfonts}
\usepackage{amsthm}
\usepackage{mathtools}
\usepackage{quiver}
\usepackage{hyperref}
\usepackage{xcolor}
\usepackage[margin=1in]{geometry}
\usepackage{comment}
\usepackage{mathdots}
\usepackage{multirow}

\usepackage[style=numeric]{biblatex}
\bibliography{references}

\newtheorem{thrm}{Theorem}[subsection]
\newtheorem{lem}[thrm]{Lemma}

\theoremstyle{definition}
\newtheorem{deff}[thrm]{Definition}

\theoremstyle{remark}
\newtheorem{rem}[thrm]{Remark}
\newtheorem{note}[thrm]{Note}

\newcommand{\orbit}{\mathrm{orb}}

\newcommand{\Cc}{\mathbb{C}}

\newcommand{\Nn}{\mathbb{N}}

\newcommand{\Pp}{\mathbb{P}}

\newcommand{\Rr}{\mathbb{R}}

\newcommand{\Zz}{\mathbb{Z}}
\title{Arithmetic Constraints on Charges of Symmetric Monopoles}

\author{Paul Cusson}

\date{October 2026}

\begin{document}

\maketitle

\begin{abstract}
    For Euclidean BPS monopoles with gauge group $SU(N)$, maximal symmetry breaking, and sufficiently generic spectral data, we produce necessary congruence relations that their magnetic charges must satisfy when $N \geq 3$. As such, a monopole as above with given charges can have certain subgroups of $SO(3)$ ruled out as possible groups of symmetries. In particular, the axially symmetric case may only admit non-trivial magnetic charges in the non-generic setting.
\end{abstract}

\tableofcontents

\section{Introduction}

The objects of interest in these notes are Euclidean BPS $SU(N)$-monopoles, with maximal symmetry breaking, and with generic spectral data. See the article \cite{hurtubise1989construction} by Hurtubise and Murray for the relevant background material. Given a closed subgroup $G$ of $SO(3)$, much work has been done on classifying those monopoles that have $G$ as their group of symmetry, although most of the results have been worked out only for $SU(2)$-monopoles; see \cite{braden2023towards}, \cite{braden2010tetrahedrally}, \cite{hitchin1983construction}, \cite{hitchin1995symmetric}. The genericity of spectral data only comes into play when $N \geq 3$, where the spectral curve decomposes into $N-1$ curves corresponding to the magnetic charges $m_1, \ldots, m_{N-1}$ so that we may consider transverse intersections. This is the origin of the strategy used in these notes, when we choose $G$ to be finite. We recall the following classical result (see Theorem 6.12.1 in \cite{artin2011algebra}).

\begin{thrm}
    A finite subgroup of $SO(3)$ is one of the following up to conjugation:
    \begin{enumerate}
        \item $\Zz_n$, the cyclic group of rotations by multiples $2\pi/n$ about a centered line,
        \item $D_n$, the dihedral group of symmetries of a regular $n$-gon, with rotations by multiples of $2\pi/n$ about a centered line, and a rotation by $\pi$ about a line joining the origin with one of the vertices of the $n$-gon,
        \item $H_{12}$, the tetrahedral group of 12 rotational symmetries of a tetrahedron,
        \item $H_{24}$, the octahedral group of 24 rotational symmetries of a cube or an octahedron,
        \item $H_{60}$, the icosahedral group of 60 rotational symmetries of a dodecahedron or an icosahedron.
    \end{enumerate}
\end{thrm}

The results in these notes all take the following form. Given a monopole with finite symmetry group $G$ and generic spectral data, its charges must satisfy arithmetic constraints of the form
$$ 2m_i m_j \equiv a \pmod b, $$
for integers $a$ depending on the orbits in the spectral curve's components and $b$ depending on the order of $G$.
\newline

\textbf{Organisation of the notes.} Section $2$ sets the notation and briefly states the relevant results on monopoles and spectral curves used here. Section $3$ presents and proves the results for all finite subgroups $G$ of $SO(3)$ up to conjugation. Section $4$ presents examples of charges ruled out for generic monopoles with given symmetry groups. Section $5$ discusses the case of $G = S^1$ and how it relates to the finite group cases.
\newline

\textbf{Acknowledgments.} These notes grew from work on my PhD thesis during the fall of 2025 and winter of 2026. I would like to thank my supervisors Professor Benoit Charbonneau and Professor Ruxandra Moraru for their guidance, suggestions, and many helpful meetings and discussions. I am also grateful to Dr. Ababacar Sadikhe Djité. My study of $SU(N)$-monopoles had kick-started while working with him. His visit at the University of Waterloo during the summer of 2025 and our near-daily discussions proved fruitful and productive in my understanding of the seminal works referenced herein.
\newline

\textbf{Use of generative AI.} No LLMs or other generative AI models were used in the work or writing of these notes.

\section{Preliminaries}

We shall denote a Euclidean BPS $SU(N)$-monopole simply as an $SU(N)$-monopole, or just as a monopole when there is no room for ambiguity. We assume all monopoles to have maximal symmetry breaking.

We identify $T\Pp^1$ with the mini-twistor space of $\Rr^3$ whose elements are oriented lines. Such an element $\gamma$ may be denoted as $\gamma(t) = tu + v$ for $(u,v) \in S^2 \times \Rr^3$ and $u \cdot v =0$, or as a pair $(\eta, \zeta) \in \Cc^2$ of local coordinates. Here $\zeta$ is in the standard chart of $\Pp^1$ and $\eta$ is the tangent space coordinate. We have an involution $\tau$ on $T\Pp^1$ that exchanges the orientation of lines $\gamma$. Together with the fact that the $SO(3)$ action does not depend on the orientation of the lines, we have that the action of $SO(3)$ and $\tau$ on $T\Pp^1$ commute.

Let $S \subset T\Pp^1$ be the spectral curve of an $SU(N)$-monopole with magnetic charges $m_1$ to $m_{N-1}$, and let $P(\eta, \zeta)$ be the polynomial cutting out $S$. We have
$$ S= S_1 \cup \cdots \cup S_{N-1}. $$
The $S_i$ are compact algebraic curves belonging to the divisor class $|\mathcal{O}_{T\Pp^1}(2m_i)|$, and they are $\tau$-invariant. Let $P_i(\eta, \zeta)$ be the factor of $P$ cutting out $S_i$. The monopole is called \textit{generic} if
$$ |S_i \cap S_{i+1}| = 2m_i m_{i+1} $$
for all $1 \leq i \leq N-2$; in particular, they intersect transversely. In this case, $S_i \cap S_{i+1} = S_{i,i+1} \cup S_{i+1, i}$ as a disjoint union into equal size parts, such that $\tau(S_{i,i+1}) = S_{i+1,i}$. We note that this definition of genericity, introduced in \cite{murray1984non}, is only defined for consecutive components of $S$ due to a correspondence with the Dynkin diagram of $SU(N)$. The stronger condition that any two components $S_i$ and $S_j$ intersect transversely is also Zariski-open in the moduli space of monopoles, as discussed in \cite{murray1984non} before Lemma 52. Thus there is still an abundance of such monopoles. We shall assume this stronger genericity condition to be the case from here on out, giving the following.

\begin{lem} \label{lemma disjoint S_i,i+1}
    Assume all components $S_i$ of the spectral curve intersect transversely. Then
    $$ |S_i \cap S_{j}| = 2m_i m_j $$
    for any $i \neq j$.
\end{lem}

\begin{proof}
    The discussion before Lemma 1.8 in \cite{hurtubise1989construction} regarding approximate solutions of $SU(N)$ monopoles via superposition of well separated charge 1 $SU(2)$-monopoles applies. The same intersection theory argument via computing degrees can then be used for $|S_i \cap S_j|$ as with $|S_i \cap S_{i+1}|$.
\end{proof}
    

\section{Finite symmetry groups}

If any non-trivial element $A \in SO(3)$ sends $S_1$ to $S_2$, say, then not only do the two curves have the same degree, $S_1$ and $S_2$ must also be isometric, under the metric on $T\Pp^1$ inherited from its inclusion in $T\Rr^3$. This is highly non-generic, so we shall also assume henceforth that the symmetry groups of monopoles fix their spectral curve components $S_i$. We compile all assumptions made on monopoles in these notes for clarity and ease of access:
\newline

\begin{tabular}{|p{15cm}|}
    \hline
    \textbf{Assumptions:}\\
    \hline
    1. Monopoles are Euclidean BPS with gauge group $SU(N)$ and maximal symmetry breaking. \\
    2. The components $S_i$ of their spectral curves intersect pair-wise transversely. \\
    3. The monopoles' symmetry groups fix the components. \\
    \hline
\end{tabular}
\newline \\

\begin{deff}
    A monopole is \textbf{nice} if it satisfies all the above assumptions.
\end{deff}

Consider a nice monopole invariant under a finite group $G$ of order $n$. Let $\gamma \in S_i \cap S_j$ and let $\orbit(\gamma)$ be the orbit of $\gamma$ under the action of $G$. By the orbit-stabilizer theorem, the size of $\orbit(\gamma)$ divides $n$. We can therefore list the possible sizes of the orbit decomposition of $S_i \cap S_j$ with respect to the action of $G$ as certain integer partitions of $n$. Given the magnetic charges of the monopole, we may rule out possible orbit decompositions, and, furthermore, certain symmetry groups entirely.

\begin{lem} \label{lemma tau invariant orbits}
    Let $\gamma \in S_i \cap S_j$. Then $\tau(\gamma) \in \orbit(\gamma)$ if and only if $\orbit(\gamma)$ is $\tau$-invariant.
\end{lem}

\begin{proof}
    Let $\tau(\gamma) \in \orbit(\gamma)$. Since $\gamma \neq \tau(\gamma)$, the orbit of $\gamma$ has at least two elements. Let $\gamma'$ be any other element in the orbit, and let $A \in G$ be such that $A \cdot \gamma = \gamma'$. Then $\tau(\gamma') = A \cdot \tau(\gamma)$, so $\tau(\gamma) \in \orbit(\gamma)$ if and only if $\tau(\gamma') \in \orbit(\gamma)$ for all $\gamma' \in \orbit(\gamma)$.
\end{proof}

We thus have two cases for the orbit decomposition of $S_i \cap S_{i+1}$. If an orbit is contained in only one of the subsets $S_{i, i+1}$ or $S_{i+1,i}$, there is a mirror copy of that orbit in the other subset given by $\tau$. Otherwise, the orbit has an even number of elements, half of which lie in each subset.

\subsection{Cyclic symmetry group}

Consider a nice monopole invariant under $\Zz_n \subset SO(3)$ acting as rotations around the $z$-axis. The only lines in $\Rr^3$ that can be fixed by this action are the two orientations of the $z$-axis. Thus there can be at most two fixed points of the action, i.e. at most two orbits of size 1. For any other oriented line $\gamma$, its orbit has $n$ elements. This gives us the following.

\begin{lem} \label{lemma z-axis cyclic}
    For $i \neq j$, the orbit decomposition of the $2m_i m_j$ points of $S_i \cap S_j$ of the monopole has either 0 or 2 fixed points, with all other orbits of size exactly $n$.
\end{lem}

We may also look at the cases where an orbit is $\tau$-invariant.

\begin{lem} \label{lemma involution in orbit}
    An element $\gamma \in S_i \cap S_j$ satisfies $\tau (\gamma) \in \orbit(\gamma)$ if and only if $n$ is even and, writing $\gamma$ as $(u,v) \in S^2 \times \Rr^3$, we have $v = (0,0,z)$ for some $z \in \Rr$.
\end{lem}

\begin{proof}
    If $\tau(\gamma) \in \orbit(\gamma)$, there exists some rotation $A \in G$ about the $z$-axis such that $A \cdot \gamma = \tau(\gamma)$. This is possible if and only if $\gamma$ intersects the $z$-axis, is parallel to the $xy$-plane, and $A$ is a rotation by $\pi$. A rotation by $\pi$ in $G$ occurs precisely when $n$ is even, and $\gamma$ is of the above form precisely when $v = (0,0,z)$.
\end{proof}

\begin{note}
    When $n$ is even and $\tau(\gamma) \in \orbit(\gamma)$, each of $S_{i,i+1}$ and $S_{i+1,i}$ contain half the elements of $\orbit(\gamma)$.
\end{note}

These impose number theoretical restrictions on the possible values of each $m_i$, giving the following result.

\begin{thrm} \label{thrm cyclic symmetry}
    A necessary condition for a nice monopole to be invariant under the subgroup $\Zz_n \subset SO(3)$ of rotations about the $z$-axis is that, for each $i \neq j$,
    \begin{equation} \label{number theory condition}
        2m_i m_j \equiv 0 , 2 \pmod{n}.
    \end{equation}
    The first case occurs generically, and the second case occurs only when the two orientations of the $z$-axis belong in $S_i \cap S_j$. 
\end{thrm}

\begin{proof}
    We have that $|S_i \cap S_j|=2m_i m_j$, and by Lemma \ref{lemma z-axis cyclic} it either contains 0 or 2 fixed points, which can only be the two orientations of the $z$-axis. All other orbits have size $n$. Thus either $2m_i m_j$ or $2m_i m_j-2$ is a multiple of $n$, as desired.

\end{proof}

\begin{rem}
    If $n$ is odd, condition \ref{number theory condition} reduces to
    $$ m_i m_j \equiv 0,1 \pmod{n}. $$
    This can also be seen from the fact that, for odd $n$, no orbit of size $n$ in $S_i \cap S_j$ is $\tau$-invariant by lemma \ref{lemma involution in orbit}, so every orbit of every size has a mirror copy.
\end{rem}

\begin{rem}
    The $z$-axis, as an element of $T\Pp^1$, corresponds to the point $(\eta, \zeta) = (0,0)$. Thus a spectral curve $S_i$ contains the $z$-axis if and only if the polynomial $P_i (\eta, \zeta)$ cutting out $S_i$ has no constant term. 
\end{rem}

\begin{rem}
    There are $\binom{N-1}{2}$ values $2m_i m_j$ to consider, as opposed to just $N-2$ originally if we only assume that $S_i, S_{i+1}$ intersect transversely for all $1 \leq i \leq N-2$. This is far stronger as $N$ increases.
\end{rem}

\subsection{Dihedral symmetry group}

In this case, we need to consider the additional reflection element of $D_n$, which is still a rotation in $\Rr^3$. We shall fix this rotation to be by $\pi$ about the $x$-axis. As a result, there are no longer any fixed points, since the $z$-axis changes orientation after a reflection.

The action of $D_n$ is given by the rotation matrices
$$ \Theta_k = \begin{pmatrix}
    \cos(\theta_k) & -\sin(\theta_k) & 0\\
    \sin(\theta_k) & \cos(\theta_k) & 0\\
    0 & 0 & 1
\end{pmatrix}, $$
where $\theta_k = \frac{2 \pi k}{n}$ for $0 \leq k \leq n-1$, and by the reflection matrix
$$ R = \begin{pmatrix}
    1 & 0 & 0\\
    0 & -1 & 0\\
    0 & 0 & -1
\end{pmatrix}. $$
Given a line $\gamma = (u,v) \in S^2 \times \Rr^3$, for any matrix $A \in D_n$ the action is $A \cdot \gamma = (Au, Av)$.

Let $v_1, \ldots, v_n$ be the vertices of a regular $n$-gon in the $xy$-plane, all unit vectors, with $v_1 = (1,0,0)$. There are $2n$ oriented lines of interest that are induced by these vertices. If $n$ is odd, for each vector $v_i$ we have two oriented lines given by the span of $v_i$. When $n$ is even, the vectors $v_1, \ldots, v_{n/2}$, and $v_i + v_{i+1}$ for $1 \leq i \leq n$ each give two oriented lines. In either case, we'll call this set of $2n$ oriented lines $\Gamma$.

\begin{lem} \label{lemma dihedral orbit sizes}
    Let $\gamma$ be an oriented line in $\Rr^3$. If $\gamma$ is either orientation of the $z$-axis, then with respect to the action of $D_n$, $\orbit(\gamma)$ has $2$ elements. If $\gamma \in \Gamma$, $\orbit(\gamma)$ has $n$ elements. In all other cases, $\orbit(\gamma)$ has $2n$ elements.
\end{lem}

\begin{proof}
    Let $\gamma = (u,v) \in S^2 \times \Rr^3$. If $\gamma$ is either orientation of the $z$-axis, $u = (\pm 1, 0,0)$ and $v=(0,0,0)$, and only the element $R \in D_n$ acts non-trivially on $\gamma$, giving 2 elements in the orbit.

    It is clear that $\Gamma$ is $D_n$-invariant, but the action is not transitive on $\Gamma$. For odd $n$, $\gamma$ and $-\gamma$ are not related, and for even $n$, the lines spanned by $v_i$ are not related to those spanned by $v_j + v_{j+1}$. This partitions $\Gamma$ into two orbits of size $n$ induced by $\Zz_n$.

    If $v \neq (0,0,0)$, $\gamma$ is neither the $z$-axis nor contained in $\Gamma$. The $n$ elements $\Theta_k$ act as $\Zz_n$ and produce $n$ distinct lines. If $v \neq (\lambda, 0,0)$, then $Rv \neq v$, and if $v = (\lambda, 0,0)$, then $Ru \neq u$. In either case there are $2n$ orbit elements.
\end{proof}

\begin{thrm} \label{thrm dihedral symmetry}
    If a nice monopole is invariant under the dihedral group $D_n \subset SO(3)$, defined about the $z$-axis, then the following holds for each pair $i \neq j$.
    \begin{enumerate}
        \item If both the $z$-axis and an element of $\Gamma$ belong to $S_i \cap S_j$, $$ 2m_i m_j \equiv 2 \pmod{n}. $$
        \item If the $z$-axis but no element of $\Gamma$ belongs to $S_i \cap S_j$, $$ 2m_i m_j \equiv 2 \pmod{2n}. $$
        \item If an element of $\Gamma$ belongs to $S_i \cap S_j$ but the $z$-axis does not, $$ 2m_i m_j \equiv 0 \pmod{n}. $$
        \item If neither the $z$-axis nor any element of $\Gamma$ belong to $S_i \cap S_j$, $$ m_i m_j \equiv 0 \pmod{n}. $$
    \end{enumerate}
    The last case occurs generically.
\end{thrm}

\begin{proof}
    The four possible cases follow from the possible orbits, and the congruences follow from lemma \ref{lemma dihedral orbit sizes}.

\end{proof}

\subsection{Tetrahedral symmetry group}

Let $H_{12} \subset SO(3)$ be the subgroup of rotational symmetries of a regular tetrahedron $T \subset \Rr^3$ whose vertices are unit vectors and contains the vector $(0,0,1)$. The group $H_{12}$ is isomorphic to $A_4 \subset S_4$, the even permutations on four elements. Labeling the vertices of $T$ from 1 to 4 makes the action of $A_4$ clear.

\begin{lem}
    Let $\gamma$ be an oriented line in $\Rr^3$. The possible sizes of the orbit of $\gamma$ under the action of $H_{12}$ are 4, 6, and 12.
\end{lem}

\begin{proof} \label{lemma tetrahedron orbit sizes}
    Of the factors of 12, we easily rule out 1, 2, and 3. A line spanned by a vertex of $T$ has 4 elements in its orbit. A line joining the origin to the midpoint of an edge of $T$ has 6 elements in its orbit. Any other line containing the origin has an orbit of size 12.
    
    If a line $\gamma = (u,v)$ does not pass through the origin, then $v \neq 0$. The orbit of $v$ falls into the same three categories as above. If it has size 4, then exactly 3 elements of $H_{12}$ fix $v$, the rotations about the line spanned by $v$, colinear with one of the vertices of $T$. The action of these rotations on $u$ must then have order 3, giving 12 elements in the orbit of $\gamma$. Similarly, when the orbit of $v$ has 6 elements, $v$ is colinear with the midpoint of one of the edges, and only the identity and the rotation by $\pi$ about about that midpoint fix $v$. This rotation does not fix $u$, giving 12 elements in the orbit of $\gamma$ once more. So in all three cases, $\orbit(\gamma)$ has size 12.
\end{proof}

We thus have three special finite orbits, $V^+$, $V^-$, and $E$, with $V^+$ corresponding to the lines spanned by the vertices of $T$. We have $V^- = \tau(V^+)$, and $E$ consists of the lines spanned by the midpoints of the edges of $T$. There are 4 elements in each of $V^+$ and $V^-$, and 6 elements of $E$. Note that $\tau(E) = E$.

\begin{thrm} \label{thrm tetrahedral symmetry}
    If a nice monopole is $H_{12}$-invariant, then the following cases are exhaustive for every pair $i \neq j$.
    \begin{enumerate}
        \item $S_i \cap S_j$ contains $V^+$, $V^-$, and $E$, and
        $$ m_i m_j \equiv 1 \pmod{6}. $$
        \item $S_i \cap S_j$ contains $V^+$ and $V^-$ but no element of $E$, and
        $$ m_i m_j \equiv 4 \pmod{6}. $$
        \item $S_i \cap S_j$ contains $E$ but no elements of either $V^+$ or $V^-$, and
        $$ m_i m_j \equiv 3 \pmod{6}. $$
        \item $S_i \cap S_j$ contains no elements of $V^+$, $V^-$, or $E$, and
        $$m_i m_j \equiv 0 \pmod{6}.$$
    \end{enumerate}
    The last case occurs generically.
\end{thrm}

\begin{proof}
    First, if any element of $V^+$ or $V^-$ belongs to $S_i \cap S_j$, then $V^+ \cup V^- \subset S_i \cap S_j$, since $S_i \cap S_j$ is $H_{12}$- and $\tau$-invariant. The same holds for elements in $E$. Thus the four cases in the theorem cover all possibilities.
    
    From the possible sizes of the orbits in the orbit decomposition of $S_i \cap S_j$, there exists non-negative integers $a,b,c$ such that
    $$ 2m_i m_j = 4a +6b +12c. $$
    But the only possibilities for $a$ are 0 and 2, and for $b$, 0 and 1, given that the only orbits of size 4 or 6 are $V^+$, $V^-$, and $E$. The congruences for each of the four cases then follow easily.
\end{proof}

\subsection{Octahedral symmetry group}

Let $H_{24} \subset SO(3)$ be the subgroup of rotational symmetries of a regular octahedron $O \subset \Rr^3$ whose vertices are the vectors $(\pm 1, 0, 0)$, $(0, \pm 1, 0)$, $(0,0, \pm 1)$. This is also the symmetry group of the cube. The group $H_{24}$ is isomorphic to $S_4$, the symmetry group on four elements.

\begin{lem}
    Let $\gamma$ be an oriented line in $\Rr^3$. The possible sizes of the orbit of $\gamma$ under the action of $H_{24}$ are 6, 8, 12, and 24.
\end{lem}

\begin{proof}
    The proof is similar to the tetrahedral case. Of the factors of 24, we can easily rule out orbits of size 1, 2, 3, and 4. A line spanned by a vertex, an edge, or a face of $T$ has 6, 12, or 8 elements in its orbit respectively. Any other line containing the origin has an orbit of size 24.

    Let $\gamma = (u,v)$ be a line that does not pass through the origin, i.e. $v \neq 0$. The orbit of the vector $v$ falls into the above four categories. Similarly as with $H_{12}$, in each case the number of elements of $H_{24}$ fixing $v$ and the corresponding orbit size of $u$ multiply to 24. This gives 24 elements in the orbit of $\gamma$.
\end{proof}

We have three special finite orbits as well, $V$, $E,$ and $F$, corresponding to the vertices, edges, and faces of $O$. In each case, they are invariant under $\tau$, as opposed to the two vertex orbits of $T$.

\begin{thrm} \label{thrm octahedral symmetry}
    Consider an $H_{24}$-invariant nice monopole. For every pair $i \neq j$, each of $V$, $E$, and $F$ is either contained in $S_i \cap S_j$ or is disjoint from it. From these 8 cases, let $\{ n_V, n_E, n_F \}$ take values in $\{ 0,1 \}$ depending on whether the corresponding orbits belong to $S_i \cap S_j$. Then we must have that
    $$ m_i m_j \equiv 3n_V + 6 n_E + 4n_F \pmod{12}. $$
    The case with $n_V = n_E = n_F =0$ occurs generically.
\end{thrm}

\begin{proof}
    It is clear that $V$, $E$, and $F$ are contained in $S_i \cap S_j$ whenever they intersect it. Given that these orbits are the only ones with sizes 6, 12, and 8, respectively, there can be at most one such orbit in $S_i \cap S_j$. This imposes
    $$2m_i m_{i+1} = 6n_V + 12n_E + 8n_F +24a$$
    for some non-negative integer $a$, which reduces to the congruence in the theorem.
\end{proof}

\subsection{Icosahedral symmetry group}

Let $H_{60} \subset SO(3)$ be the subgroup of rotational symmetries of a regular icosahedron $I \subset \Rr^3$. The icosahedron is embedded such that its vertices are unit length, they include the vectors $(0,0,\pm 1)$, and its intersection with the $xy$-plane is a regular octagon, oriented so as to have a vertex on the $x$-axis. $H_{60}$ is also the symmetry group of the dodecahedron, and is isomorphic to $A_5$, the alternating subgroup of even permutations in $S_5$.

\begin{lem}
    Let $\gamma$ be an oriented line in $\Rr^3$. The possible sizes of the orbit of $\gamma$ under the action of $H_{60}$ are 12, 20, 30, and 60.
\end{lem}

\begin{proof}
    Again, we rule out the small factors of 60, and the minimal orbit size is 12, achieved by the lines spanned by the vertices of $I$. Because $H_{60}$ contains rotations by $\pi$, any orbit must have an even number of elements. This rules out 15. The lines spanned by the center of the faces of $I$ have orbit size 20, and those spanned by the midpoints of the edges have orbit size 30. Any other line containing the origin has 60 elements in its orbit.

    Let $\gamma = (u,v)$ be a line that does not pass through the origin, so that $v \neq 0$. The same case analysis as the tetrahedral and octahedral cases applies, and they all have orbits of size 60.
\end{proof}

Like before, this gives us three special finite orbits, $V$, $E$, and $F$, corresponding to the vertices, edges, and faces of $I$, all $\tau$-invariant.

\begin{thrm} \label{thrm icosahedral symmetry}
    Consider a $H_{60}$-invariant nice monopole. For every pair $i \neq j$, each of $V$, $E$, and $F$ is either contained in $S_i \cap S_j$ or is disjoint from it. From these 8 cases, let $\{ n_V, n_E, n_F \}$ take values in $\{0,1\}$ depending on whether the corresponding orbits belong to $S_i \cap S_j$. Then
    $$ m_i m_j \equiv 6n_V +15n_E +10n_F \pmod{30}. $$
    The case with $n_V = n_E = n_F =0$ occurs generically.
\end{thrm}

\begin{proof}
    The congruence in this case follows from
    $$ 2m_i m_{i+1} = 12n_V + 30n_E + 20n_F + 60a. $$ The rest of the proof follows entirely that of the octahedral case.
\end{proof}

\section{Constrained charges in the generic cases}

We consider, in this section, nice $SU(4)$-monopoles and the charges $(m_1, m_2, m_3)$ that each symmetry group's \textit{generic} congruence condition force. We begin with the symmetry groups of the Platonic solids.
\newline

\begin{tabular}{ | p{3cm} | p{6cm} | p{6cm} | }
    \hline
    Symmetry group & generic congruence condition for all $i \neq j$ & possible magnetic charges up to permutation, with $a,b,c \in \Nn$ \\
    \hline\hline
    Tetrahedral & $m_i m_j \equiv 0 \pmod 6$ & $(a, 6b, 6c)$ \\
    & & $(2a, 3b, 6c)$\\
    \hline
    Octahedral & $m_i m_j \equiv 0 \pmod {12}$ & $(a, 12b, 12c)$ \\
     & & $(2a, 6b, 12c)$ \\
     & & $(3a, 4b, 12c)$ \\
    \hline
    Icosahedral & $m_i m_j \equiv \pmod {30}$ & $(a, 30b, 30c)$ \\
     & & $(2a, 15b, 30c)$ \\
     & & $(3a, 10b, 30c)$ \\
     & & $(5a, 6b, 30c)$ \\
    \hline
\end{tabular}

For $G = \Zz_n$ as in Theorem \ref{thrm cyclic symmetry} with fixed $n$, the generic congruence condition is
$$ 2 m_i m_j \equiv 0 \pmod n $$
for all $i \neq j$. For $G = D_n$ as in Theorem \ref{thrm dihedral symmetry}, the condition is
$$ m_i m_j \equiv 0 \pmod n. $$
For composite $n$, the possible magnetic charges have few restrictions. On the other hand, when $n$ is a prime $p \neq 2$ for $\Zz_p$, and for any prime $p$ for $D_p$, the set of possible $(m_1, m_2, m_3)$ are up to permutation of the form
$$ (a, bp, cp) $$
for any integers $a, b, c \geq 1$.

\section{Axially symmetric monopoles}

When $G = S^1$ acts about the $z$-axis, the intersection of any two components $S_i$ and $S_j$ are the points $(0,0)$ and $(0,\infty)$. But these points correspond precisely to the two orientations of the $z$-axis. Since $\Zz_n \subset S^1$, we apply the \textit{non-generic} congruence condition
$$ 2 m_i m_j \equiv 2 \pmod n $$
for every positive integer $n$.
This forces all $m_i =1$ when the monopole satisfies our initial assumptions of genericity. In particular, the only examples of axially symmetric, maximal symmetry breaking, $SU(N)$-monopoles with non-trivial magnetic charge cannot be generic, even merely in the sense of \cite{murray1984non}. The moduli space of such axially symmetric monopoles must then lie in the complement of the Zariski-open subspace of generic monopoles.

\nocite{*}
\printbibliography

\end{document}